\documentclass[11pt]{amsart}
\usepackage{amsmath,amsfonts,amsthm,amsopn,amssymb}
\usepackage{cite,marginnote}
\usepackage{graphicx}
\usepackage{graphics}
\usepackage{subfigure}

\usepackage{color,enumitem,graphicx}
\usepackage[colorlinks=true,urlcolor=blue,
citecolor=red,linkcolor=blue,linktocpage,pdfpagelabels,
bookmarksnumbered,bookmarksopen]{hyperref}
\usepackage[english]{babel}

\usepackage[left=2.7cm,right=2.7cm,top=2.8cm,bottom=2.8cm]{geometry}
\usepackage[hyperpageref]{backref}

\numberwithin{equation}{section}

\makeindex
\newtheorem{theorem}{Theorem}[section]
\newtheorem{proposition}[theorem]{Proposition}
\newtheorem{lemma}[theorem]{Lemma}
\newtheorem{remark}[theorem]{Remark}
\newtheorem{example}[theorem]{Example}
\newtheorem{corollary}[theorem]{Corollary}
\newtheorem{definition}[theorem]{Definition}

\newcommand{\ud}{\mathrm{d}}
\newcommand{\RN}{\mathbb R^N}

\newcommand{\iy}{\infty}

\newcommand{\s}{\section}

\newcommand{\DD}{\Delta}
\newcommand{\g}{\gamma}
\newcommand{\G}{\Gamma}

\newcommand{\R}{\mathbb R}
\newcommand{\al}{\alpha}

\newcommand{\ti}{\tilde}

\newcommand{\re}[1]{(\ref{#1})}

\newcommand{\rg}{\rightarrow}

\newcommand{\lab}{\label}
\newcommand{\bt}{\begin{theorem}}
\newcommand{\et}{\end{theorem}}
\newcommand{\bl}{\begin{lemma}}
\newcommand{\el}{\end{lemma}}
\newcommand{\bd}{\begin{definition}}
\newcommand{\ed}{\end{definition}}
\newcommand{\bc}{\begin{corollary}}
\newcommand{\ec}{\end{corollary}}
\newcommand{\bp}{\begin{proof}}
\newcommand{\ep}{\end{proof}}
\newcommand{\bx}{\begin{example}}
\newcommand{\ex}{\end{example}}
\newcommand{\bi}{\begin{exercise}}
\newcommand{\ei}{\end{exercise}}
\newcommand{\bo}{\begin{proposition}}
\newcommand{\eo}{\end{proposition}}
\newcommand{\br}{\begin{remark}}
\newcommand{\er}{\end{remark}}
\newcommand{\be}{\begin{equation}}
\newcommand{\ee}{\end{equation}}
\newcommand{\ba}{\begin{align}}
\newcommand{\ea}{\end{align}}
\newcommand{\bn}{\begin{enumerate}}
\newcommand{\en}{\end{enumerate}}
\newcommand{\bg}{\begin{align*}}
\newcommand{\bcs}{\begin{cases}}
\newcommand{\ecs}{\end{cases}}

\newcommand{\bean}{\begin{eqnarray*}}
\newcommand{\eean}{\end{eqnarray*}}

\title[Fractional Kirchhoff equation]{Fractional Kirchhoff equation with a general critical nonlinearity}

\author[H.\ Jin]{Hua Jin}
\author[W. B.\ Liu]{Wenbin Liu}

\address[H.\ Jin]{\newline\indent College of Science
\newline\indent
China University of Mining and Technology
\newline\indent
Xuzhou, 221116, China}
\email{\href{mailto:huajin@cumt.edu.cn}{huajin@cumt.edu.cn}}

\address[W. B.\ Liu]{\newline\indent College of Science
\newline\indent
China University of Mining and Technology
\newline\indent
Xuzhou, 221116, China}
\email{\href{mailto:liuwenbin-xz@163.com}{liuwenbin-xz@163.com}}

\thanks{W. Liu is the corresponding author. }


\begin{document}

\begin{abstract}
In this paper, we study the fractional Kirchhoff equation with critical nonlinearity
\begin{align*}
\left(a+b\int_{\RN}|(-\Delta)^{\frac{s}{2}}u|^2dx\right)(-\Delta)^su+u=f(u)\ \   \mbox{in}\ \ \RN,
\end{align*}
where  $N>2s$ and $(-\Delta)^s$ is the fractional Laplacian with $0<s<1$. By using a perturbation approach, we prove the existence of solutions to the above problem  without the Ambrosetti-Rabinowitz condition when the parameter $b$ small. What's more, we obtain the asymptotic behavior of solutions as $b\rg 0$.
\end{abstract}
\maketitle

\s{Introduction and main result}
In this paper, we are concerned with the following fractional Kirchhoff equation
\be
\lab{fractional kirchhoff}
\left(a+b\int_{\RN}|(-\Delta)^{\frac{s}{2}}u|^2dx\right)(-\Delta)^su+u=f(u)\ \   \mbox{in}\ \ \RN,
\ee
where $N>2s$ with $0<s<1$, $a,b$ are positive constants and $(-\Delta)^su$ is the fractional Laplacian which arises in the description of various phenomena in the applied science, such as the phase transition \cite{Y.Sire}, Markov processes \cite{App} and fractional quantum mechanics \cite{Lask2}.
When $a=1$ and $b=0$, (\ref{fractional kirchhoff}) becomes the fractional Schr\"{o}dinger equations which have been studied by many authors.
%\be\lab{fracsch}
%(-\Delta)^su+u=f(u)\ \   \mbox{in}\ \ \RN,
%\ee
%For such a class of fractional problems, many results on the existence were obtained. In \cite{Changxiaojun1}, X. Chang and Z.-Q. Wang proved the existence of ground states in the subcritical case.  In the critical case, J. Zhang et al. \cite{zjmf} showed the existence of ground states when $f$ satisfies the Berestycki-Lions type conditions. For the more general fractional Laplacian equation
%$$
%(-\Delta)^su+V(x)u=f(x,u)\ \   \mbox{in}\ \ \RN,
%$$
%there has been a lot of results,
We refer the readers to \cite{Alves,Autuori,c.brandle,Changxiaojun1} and the references therein for the details. When $s=1$, the problem (\ref{fractional kirchhoff}) reduces to the well-known Kirchhoff equation
\be
\lab{kirchhoff}
-\left(a+b\int_{\RN}|\nabla u|^2dx\right)\Delta u+u=f(u)\ \   \mbox{in}\ \ \RN,
\ee
which has been studied in the last decade,  see \cite{Nyamoradi,dyps,X.H.W.Z1}. The equation (\ref{kirchhoff}) is related to the stationary analogue of the Kirchhoff equation
$
u_{tt}-\left(a+b\int_\Omega|\nabla u|^2dx\right)\Delta u=f(x,u)
$
on $\Omega\subset \RN$  bounded, which was proposed by Kirchhoff \cite{Kirchhoff} in 1883 as a generalization the classic D'Alembert's wave equation
%$$
%\rho u_{tt}-\left(\frac{P_0}h+\frac{E}{2L}\int_0^L|u_x|^2dx\right)u_{xx}^2=f(x,u)
%$$
for free vibrations of elastic strings.
 %where $L$ is the length of the string, $h$ is the area of the cross section, $E$ is the Young modulus of the material, $\rho$ is the mass density and $P_0$ is the initial tension.
%Kirchhoff's model takes into account the changes in length of the string produced by transverse vibrations.
%Besides, (\ref{kirchhoff}) can models several physical systems\cite{Alveskirc}, where $u$ describe a process which depends on the average of itself.
%In \cite{M.chipot}, the application of (\ref{kirchhoff}) is found in biological systems, describing the growth and movement of a particular species. %The movement, modeled by the integral term, is assumed to be dependent on the energy of the entire system with $u$ being its population density.
%The early classic studies on Kirchhoff equations  were given by Bernstein \cite {Bernstein} and Pohoz\v{a}ev\cite{Pohokirc}.
%However, after Lions \cite{Lionskirc} applied the functional analysis approach to investigate it, Kirchhoff equations gained more and more attention.
%For more results on (\ref{kirchhoff}), we refer the readers to \cite{Nyamoradi,dyps,X.H.W.Z1} and the references therein.

Recently, in bounded regular domains of $\RN$, Fiscella and Valdinoci \cite{Fiscella} proposed the following fractional stationary Kirchhoff equation
\begin{align}
\lab{generalfrackirchhoff}
\left\{
\begin{array}{ll}
M\left(\int_{\RN}|(-\Delta)^{\frac{s}{2}}u|^2\right)(-\Delta)^su=f(x,u), \ \   \mbox{in}\ \ \Omega,\\
u=0 \ \ \ \ \mbox{in}\ \ \RN\backslash\Omega,
\end{array}
\right.
\end{align}
which models nonlocal aspects of the tension arising from nonlocal measurements of the fractional length of the string. Also in  bounded domains, Autuori et al. \cite{AutuoriFisPucci} dealt with the existence and the asymptotic behavior of non-negative solutions of a class of fractional stationary Kirchhoff equation. In the whole of $\RN$, Pucci et al. \cite{Pucciandsaldi} established the  existence and multiplicity of nontrivial non-negative entire  solutions of a stationary Kirchhoff eigenvalue problem. In the subcritical case, by using minimax arguments, Ambrosio  et al. \cite{Ambrosio} obtained the multiplicity results for (\ref{kirchhoff}) in $H_r^s(\RN)$ with $b$ small. Also in the subcritical case, without the (AR)-condition, the authors\cite{xiangzhangyang} investigated the existence of radial solutions by using the variational methods combined with a cut-off function technique. More recently, without the (AR)-condition and monotonicity  assumptions, in low dimension($N=2,N=3$), Z. Liu et al.\cite{liuzhisumarco} studied the existence of ground states in the critical case. To the best of our knowledge, there are few papers on the fractional Kirchhoff equations involving the critical growth in $\RN$ with $N>3$, because of the tough difficulties brought by the nonlocal term and the lack of compactness of the Sobolev embedding $H_r^s(\RN)\rg L^{2_s^*}(\RN)$.

Motivated by the works above, we investigate the existence of the positive solutions of \re{fractional kirchhoff} in  $\RN(N>2s)$ with the critical growth. Precisely, $f$ satisfies the following  conditions:
\begin{itemize}

\item [$(f_1)$] $f\in C^1(\R^+,\R)$, $\lim_{t\to 0}f(t)/t=0$ and $f(t)\equiv 0$ for $t\leq 0$,

\item [$({f_2})$]  $\lim_{t\to \iy}{f(t)/t^{2_s^*-1}}=1$, where $2_s^*=\frac{2N}{N-2s}$,

\item [$({f_3})$]  there exist $D>0$ and $p<2_s^*$ such that $f(t)\geq t^{2_s^*-1}+Dt^{p-1}$ for $t\geq 0$.

\end{itemize}
Our main result can read as
\bt\lab{Theorem 1} Suppose that $f$ satisfies $(f_1)-(f_3)$ and $\max\{2,2_s^*-2\}<p<2_s^*$, then  for  $b$ small, (\ref{fractional kirchhoff}) admits a nontrivial positive radial solution  $u_b$. What's more, along a subsequence, $u_b$ converges to $u$ in
$H_r^s(\RN)$ as $b\rg 0$, where $u$ is
a radial ground state to the limit problem
\be\lab{lb}
a(-\Delta)^s u+u=f(u),\ \ \ u\in H^s(\RN).
\ee
%\%begin{itemize}
%\item [(i)] there exists $b_0>0$ such that, for every $b\in (0,b_0)$,
%\re{fractional kirchhoff} admits a nontrivial positive radial solution  $u_b$,
%\item [(ii)] along a subsequence, $u_b$ converges to $u$ in
%$H_r^s(\RN)$ as $b\rg 0$, where $u$ is
%a radial ground state solution to the limit problem
%\be\lab{lb}
%a(-\Delta)^s u+u=f(u),\ \ \ u\in H^s(\RN).
%\ee
%\end{itemize}
\et

Because of the presence of the Kirchhoff term, in high dimension $N>4s$, for the energy functional $I_b(u)$ (see section 2), one has $I_b(tu)\rg +\infty$ as $t\rg +\infty$ for each $u\neq 0$. That means  Mountain pass geometry may not holds and  Mountain pass theorem may not be appropriate. To overcome this difficulty, we  use the variational method combined with the perturbation approach\cite{zjjjmms,zjmf} to get a special bounded (PS)-sequence.
On the other hand, because of the presence of the Kirchhoff term, for the bounded (PS)-sequence $\{u_n\}$, even $u_n\rg u_0$ weakly, it doesn't hold in general that $u_0$ is the critical point of the energy functional, which brings us more tough to get the compactness. We use the properties of the special (PS)-sequence and some results of the limit problem (\ref{lb}) to  recover the compactness. Moreover,we obtain the asymptotic behavior of the solutions of (\ref{fractional kirchhoff}) as $b\rg 0$.

The paper is organized as follows. Some preliminaries are presented in Section 2.  In Section 3, we construct the min-max level. In Section 4, we complete the proof of Theorem~\ref{Theorem 1}.
%In Section 2, we introduce the functional framework and some preliminary results.
%In Section 3, we construct the min-max level.
%In Section 4, we use a perturbation argument to complete the proof of Theorem~\ref{Theorem 1}.

%%%%%%%%%%%%%%%%%%%%%%%%%%%%%%%%%%%%%%%%%%%%%%%%%%%%%%%%%%%%%%%%%%%%%%%%%%%%%%%%%%%%%%%%%%%%%%%%

\s{Preliminaries and functional setting}

\renewcommand{\theequation}{2.\arabic{equation}}

\subsection{Fractional order Sobolev spaces}
\noindent The fractional Laplacian $(-\Delta)^s$ with $s\in(0,1)$ of a function $\phi:\RN\rightarrow\R$ is defined by
$
\mathcal{F}((-\Delta)^s\phi)(\xi)=|\xi|^{2s}\mathcal{F}(\phi)(\xi),
$
where $\mathcal{F}$ is the Fourier transform.
If $\phi$ is smooth enough, it can be computed by the following singular integral
$$
(-\Delta)^s\phi(x)=c_s\, \mbox{P.V.}\int_{\RN}\frac{\phi(x)-\phi(y)}{|x-y|^{N+2s}}\, \ud y,\, \ x\in\RN,
$$
where $c_s$ is a normalization constant and P.V. stands the principal value. For any $s\in(0,1)$, we consider the fractional order Sobolev space
$$
H^s(\RN)=\left\{u\in L^2(\RN): \int_{\RN}|\xi|^{2s}|\hat{u}|^2\, \ud\xi<\iy\right\},
$$
endowed with the norm $\|u\|=\left(\int_{\RN}(1+a|\xi|^{2s})|\hat{u}|^2\, \ud\xi\right)^{1/2}.$
$H_r^s(\RN)$ denotes the space of radial functions in $H^s(\RN)$, i.e. $H_r^s(\RN)=\{u\in H^s(\RN): u(x)=u(|x|)\}.$ The homogeneous Sobolev space $\mathcal{D}^{s,2}(\RN)$ is defined by
$$
\mathcal{D}^{s,2}(\RN)=\{u\in L^{2^\ast_s}(\RN): |\xi|^s\hat{u}\in L^2(\RN)\},
$$
which is the completion of $C_0^\iy(\RN)$ under the norm
$
\|u\|_{\mathcal{D}^{s,2}}^2=\|(-\Delta)^{s/2}u\|_2^2=\int_{\RN}|\xi|^{2s}|\hat{u}|^2\, \ud\xi.
$

For the further introduction on the fractional order Sobolev space, we refer to \cite{DPV}. Now, we introduce the following Sobolev embedding theorems.
%\bl[see \cite{Lions1}]\lab{l1} For $s\in(0,1)$,
%$H^s(\RN)$ is continuously embedded into $L^q(\RN)$ for $q\in[2,2^\ast_s]$ and compactly embedded into $L^q_{loc}(\RN)$ for $q\in[1,2^\ast_s)$. Moreover, $H_r^s(\RN)$ is compactly embedded into $L^q(\RN)$ for $q\in(2,2^\ast_s)$.
%\el

\bl[see \cite{Lions1,Cots,DPV}]\lab{l2} For any $s\in(0,1)$, $H^s(\RN)$ is continuously embedded into $L^q(\RN)$ for $q\in[2,2^\ast_s]$ and compactly embedded into $L^q_{loc}(\RN)$ for $q\in[1,2^\ast_s)$. $H_r^s(\RN)$ is compactly embedded into $L^q(\RN)$ for $q\in(2,2^\ast_s)$ and
$\mathcal{D}^{s,2}(\RN)$ is continuously embedded into $L^{2^\ast_s}(\RN)$, i.e., there exists $S_s>0$ such that
$
S_s\left(\int_{\RN}|u|^{2^\ast_s}\, \ud x\right)^{2/2^\ast_s}\le\int_{\RN}|(-\DD)^{\frac{s}2}u|^2\, \ud x.
$
\el
\subsection{The variational setting}

We define the energy functional $I_b:H^s(\RN)\to\R$ by
$$
I_b(u)=\frac{1}{2}\int_{\RN}\left(a|(-\DD)^{\frac{s}2}u|^2+u^2\right)\, \ud x+\frac{b}{4}\left(\int_{\RN}|(-\DD)^{\frac{s}2}u|^2\right)^2\, \ud x-\int_{\RN}F(u)\, \ud x,
$$
with $F(t)=\int_0^t f(\zeta)\, \ud \zeta$. It is standard to show that $I_b$ is of class $C^1$.
\bd\noindent
We call $u\in H^s(\RN)$ a weak solution of \re{fractional kirchhoff} if for any $\phi\in H^s(\RN)$,
$$
\left(a+b\|u\|_{\mathcal{D}^{s,2}}^2\right)\int_{\RN}(-\DD)^{\frac{s}2}u(-\DD)^{\frac{s}2}\phi\, \ud x+\int_{\RN}u\phi\, \ud x=\int_{\RN}f(u)\phi\, \ud x.
$$
\ed
\noindent Obviously, the critical points of $I_b$ are the weak solutions of \re{fractional kirchhoff}.

Similar to the proof of Brezis-Lieb Lemma in \cite{zjjjmms}, we can give the following lemma.

\bl\lab{Brezis-Lieb lemma} For $s\in (0,1)$, assume $(f_1)-(f_2)$ hold. Let $\{u_n\}\subset H^s(\RN)$ such that $u_n\to u$ weakly in $H^s(\RN)$ and a.e. in $\RN$ as $n\to \iy$, then
$
\int_{\RN}F(u_n)\to\int_{\RN}F(u_n-u)+\int_{\RN}F(u).
$
\el

When $b=0$, problem (\ref{fractional kirchhoff}) becomes the limit problem (\ref{lb}) which plays a crucial role in our paper. The energy functional of (\ref{lb}) is defined as
$$
L(u)=\frac{1}{2}\int_{\RN}\left(a|(-\DD)^{\frac{s}2}u|^2+u^2\right)\, \ud x-\int_{\RN}F(u)\, \, \ud x,\ \ u\in H^s(\RN).
$$

With the same assumptions on $f$ in Theorem \ref{Theorem 1}, it is not difficult to check that $L(u)$ satisfies the Mountain pass geometry. The Mountain pass value denoted by $c$ is defined by
$$
c=\inf_{\g\in\G_L}\max_{t\in[0,1]}L(\g(t))>0,$$
where $\G_L=\{\g\in C([0,1],H^s(\RN)),\g(0)=0,L(\g(1))< 0\}.$ In the following, we present some results of the ground states of (\ref{lb}) and the proof is similar as that in \cite{zjmf}.
%%%%%
\bo\lab{groundstate of limit}
Suppose $f$ satisfies $(f_1)-(f_3)$ and $\max\{2,2_s^*-2\}<p<2_s^*$. Let $S_r$ be the set of positive radial ground states of (\ref{lb}), then

(i) $S_r$ is not empty and $S_r$ is compact in $H_r^s(\RN)$,

(ii) $c<\frac{s}{N}(aS_s)^{\frac{N}{2s}}$ and $c$ agrees with the least energy level denoted by $E$, that is, there exists $\g\in\G_L$ such that $u\in\g(t)$ and $\max_{[0,1]}L(\g(t))=E$, where $u\in S_r$,

(iii)$u\in S_r$ satisfies  the Pohoz\v{a}ev identity
\be\label{pohozaev2}
\frac{N-2s}{2}\int_{\RN}a|(-\DD)^{s/2}u|^2\, \ud x+\frac{N}2\int_{\RN}u^2\, \ud x=N\int_{\RN}F(u)\, \ud x.
\ee
\eo
%%%%%%

%Let $S_r$ be the set of positive radial ground state solutions of (\ref{lb}), then we have the following compactness result which plays a crucial role in the proof of Theorem \ref{Theorem 1}.
%\bo\lab{compactness of ground state}
%Under the assumptions of  Theorem \ref{Theorem 1},  $S_r$ is compactness.
%\eo

\vskip0.1in

%%%%%%%%%%%%%%%%%%%%%%%%%%%%%%%%%%%%%%%%%%%%%%%%%%%%%%%%%%%%%%%%%%%%%%%%%%%%%%%%%%%%%%%%%%%%%%%%%
%%%%%%%%%%%%%%%%%%%%%%%%%%%%%%%%%%%%%%%%%%%%%%%%%%%%%%%%%%%%%%%%%%%%%%%%%%%%%%%%%%%%%%%%%%%%%%%%%
%%%%%%%%%%%%%%%%%%%%%%%%%%%%%%%%%%%%%%%%%%%%%%%%%%%%%%%%%%%%%%%%%%%%%%%%%%%%%%%%%%%%%%%%%%%%%%%%%
%%%%%%%%%%%%%%%%%%%%%%%%%%%%%%%%%%%%%%%%%%%%%%%%%%%%%%%%%%%%%%%%%%%%%%%%%%%%%%%%%%%%%%%%%%%%%%%%%
%%%%%%%%%%%%%%%%%%%%%%%%%%%%%%%%%%%%%%%%%%%%%%%%%%%%%%%%%%%%%%%%%%%%%%%%%%%%%%%%%%%%%%%%%%%%%%%%%
%%%%%%%%%%%%%%%%%%%%%%%%%%%%%%%%%%%%%%%%%%%%%%%%%%%%%%%%%%%%%%%%%%%%%%%%%%%%%%%%%%%%%%%%%%%%%%%%%

\s{The minimax level}
\renewcommand{\theequation}{3.\arabic{equation}}

In order to get a bounded (PS)-sequence by the local deformation argument, a different min-max level is needed.
\noindent Take $U\in S_r$ be arbitrary but fixed.  By the definition of $\hat{U}=\mathcal{F}(U)$, for $U_\tau(x)=U\big(\frac{x}{\tau}\big), \tau>0$, we have $\hat{U}_\tau(\cdot)=\tau^N\hat{U}(\tau\cdot)$. Thus
$
\int_{\RN}|(-\DD)^{s/2}U_\tau|^2\, \ud x=\tau^{N-2s}\int_{\RN}|(-\DD)^{s/2}U|^2\, \ud x.
$
From the Pohoz\v{a}ev identity (\ref{pohozaev2}), we obtain
$$
L(U_\tau)=\Big(\frac{a\tau^{N-2s}}{2}-\frac{N-2s}{2N}\tau^N\Big)\int_{\RN}|(-\DD)^{s/2}U|^2.
$$
So, there exists $\tau_0>1$ such that $L(U_{\tau})<-2$ for $\tau\ge \tau_0$. Set
$
D_b\equiv\max_{\tau\in[0,\tau_0]}I_b(U_\tau).
$
Noting that $I_b(U_\tau)=L(U_{\tau})+\frac{b}4\|U_{\tau}\|_{\mathcal{D}^{s,2}}^4$ and $\max_{\tau\in[0,\tau_0]}L(U_\tau)=E$, we have $D_b\rg E$ as $b\rg 0^+$.
%
%In order to define the uniformly bounded set of  the mountain paths(see below), we first present the following lemma.

\bl\lab{l4.1} There exist $b_1>0$ and  $\mathcal{C}_0>0$, such that for
any $0<b<b_1$ there hold
\begin{equation*}
I_b(U_{\tau_0})<-2,\qquad
\|U_\tau\|\le \mathcal{C}_0, \,\,\, \forall \tau\in (0,\tau_0],\qquad
\|u\|\le \mathcal{C}_0, \,\,\, \forall u\in S_r.
\end{equation*}
\el
\bp
Since $S_r$ is compact, it is easy to verify that there exists $\mathcal{C}_0>0$ such that the second and third part of the  assertion hold. It follows from
$
I_b(U_{\tau_0})\leq L(U_{\tau_0})+\frac{b}4\mathcal{C}_0^4
$
and $L(U_{\tau_0})<-2$ that the first part holds for any $0<b<b_1$, where $b_1>0$ small. The proof is completed.
\ep

\noindent Now, for any $b\in (0,b_1)$, we define a min-max value $C_b:=\inf_{\g\in \Upsilon_b}\max_{\tau\in [0,\tau_0]}I_b(\g(\tau))$,
where
\begin{align*}
\Upsilon_b=\big\{\g\in C([0,\tau_0],H_r^s(\RN)): \g(0)=0,
\g(\tau_0)=U_{\tau_0},\|\g(\tau)\|\le \mathcal{C}_0+1,\tau\in[0,\tau_0]\big\}.
\end{align*}

\bo\lab{bo4} $\lim\limits_{b\rg 0^+}C_b=E.$ \eo
\bp
For $\tau>0$, by $\|U_\tau\|^2=a\tau^{N-2s}\|U\|_{\mathcal{D}^{s,2}}^2+\tau^N\|U\|_2^2,$
 we can define $U_0\equiv0$. So $U_\tau\in \Upsilon_b$. Moreover,
$\limsup_{b\rg0^+}C_b\le\lim_{b\rg0^+}D_b=E=c.$ On the other hand,
for any $\g\in \Upsilon_b$, it follows from $L(U_{\tau_0})<-2$ that $\ti{\g}(\cdot)=\g(\tau_0\cdot)\in \G_L$. Thus, from the definition of $c$ and $C_b$, we obtain $C_b\geq E$ for any $b\in (0,b_1)$. The proof is completed.
\ep

\s{The Proof of Theorem \ref{Theorem 1}}
\renewcommand{\theequation}{4.\arabic{equation}}

For $\al,d>0$, we define
$$I_b^\al:=\{u\in H_r^s(\RN): I_b(u)\le\al\}
$$
and
$$
S^d:=\left\{u\in H_r^s(\RN):  \inf_{v\in
S_r}\|u-v\|\le d\right\}.
$$

%The following proposition is crucial to obtain a suitable
%(PS)-sequence for $I_b$ and plays a key role in our proof.

\bo\lab{bo5} Let $\{b_n\}_{n=1}^\infty$ be such that $\lim_{n\rg
\infty}b_n=0$ and $\{u_{b_n}\}\subset S^d$
with
$$
\lim_{n\rg\infty}I_{b_n}(u_{b_n})\le E\ \mbox{and}\
\lim_{n\rg\infty}I_{b_n}'(u_{b_n})=0.
$$
Then for $d$ small, there is $u_0\in S_r$, up to
a subsequence, such that $u_{b_n}\rg u_0$ strongly in $H_r^s(\RN)$.
\eo
\bp
For convenience, we write $u_{b_n}$ for $u_b$. Since $u_b\in S^d$, there exists $\ti u_b\in S_r$ such that $\|u_b-\ti u_b\|\le d$. Let $v_b=u_b-\ti u_b$. By the fact that $S_r$ is compact and $\|v_b\|\le d$, up to a subsequence, there exist $\ti u_0\in S_r$ and $v_0\in H^s(\RN)$, such that $\ti u_b\rg \ti u_0$ strongly in $H_r^s(\RN)$ and $v_b\rg v_0$ weakly in $H^s(\RN)$. Denoting $u_0=\ti u_0+v_0$, then $u_0\in S^d$ and $u_b\rg u_0$ weakly in $H_r^s(\RN)$. Next, we show $u_b\rg u_0$ strongly in $H_r^s(\RN)$.
Since $\lim_{n\rg\infty}I_{b}'(u_{b})=0$, then for any $\phi \in C_0^\infty(\RN)$,
$$
I_{b}'(u_{b})\phi=L'(u_b)\phi+b\|u_b\|_{\mathcal{D}^{s,2}}^2\int_{\RN}(-\DD)^{s/2}u_b(-\DD)^{s/2}\phi.
$$
It follows from Lemma \ref{l2} and  $u_b\in S^d$  that $L'(u_0)=0$ as $b\to 0$. Obviously $u_0\not\equiv 0$ by $u_0\in S^d$ with $d$ small. Thus $L(u_0)\ge E$.
Meanwhile, from Lemma \ref{Brezis-Lieb lemma},
$
I_b(u_b)=L(u_b)+\frac{b}4\|u_b\|_{\mathcal{D}^{s,2}}^4=L(u_0)+L(u_b-u_0)+o(1).
$
Together with $\lim_{n\rg\infty}I_{b_n}(u_{b_n})\le E$, we obtain $L(u_b-u_0)\le o(1)$. Thus, by $(f_1)-(f_2)$ and Sobolev embedding theorem, there exists constant $c_1>0$ such that $\|u_b-u_0\|^2\le c_1\|u_b-u_0\|^{2_s^*}$. If $\|u_b-u_0\|\nrightarrow 0$ as $b\rg 0$, there exists constant $c_2>0$ such that $\|u_b-u_0\|\ge c_2$ for $b$ small. On the other hand, from $\ti u_0\in S_r$ and $u_0\in S^d$, we get $\|\ti u_0-u_0\|\le d$. Then
$
\|u_b-u_0\|\le \|u_b-\ti u_b\|+\|\ti u_b- \ti u_0\|+\|\ti u_0-u_0\|\le 2d+o(1),
$
which is a contradiction for $d$ small. The proof is completed.
\ep

\br \lab{remark1}By Proposition \ref{bo5}, for small $d\in(0,1)$, there exist $\omega>0, b_0>0$ such that
\be\lab{deform}
\hbox{$\|I_b'(u)\|\ge\omega$ for
$u\in I_b^{D_b}\bigcap(S^d\setminus S^{\frac{d}{2}})$
and $b\in(0,b_0)$.}
\ee
\er
\noindent Thus, we have the following proposition.
\bo\lab{bo6} There exists $\al>0$ such that for small $b>0$ and $\g(\tau)=U(\frac{\cdot}{\tau}), \tau\in (0,\tau_0]$,
$$
I_b(\g(\tau))\ge C_b-\al\ \ \mbox{implies that}\ \ \g(\tau)\in
S^{\frac{d}{2}},
$$
\eo

\bp By the Pohoz\v{a}ev identity (\ref{pohozaev2}), $I_b(\g(\tau))=\Big(\frac{a\tau^{N-2s}}{2}-\frac{N-2s}{2N}\tau^N\Big)\|U\|_{\mathcal{D}^{s,2}}^2+\frac{b}{4}\tau^{2N-4s}\|U\|_{\mathcal{D}^{s,2}}^4.$
Then $\lim_{b\rg0^+}\max_{\tau\in[0,\tau_0]}I_b(\g(\tau))=\max_{\tau\in[0,\tau_0]}\Big(\frac{a\tau^{N-2s}}{2}-\frac{N-2s}{2N}\tau^N\Big)\|U\|_{\mathcal{D}^{s,2}}^2=E.$
On the other hand, $\lim\limits_{b\rg 0^+}C_b=E.$ The conclusion follows. \ep

\noindent
Thanks to (\ref{deform}) and Proposition \ref{bo6}, we can prove the following proposition, which assures the existence of a bounded
(PS)-sequence for $I_b$. The proof is similar as that in \cite{zjjjmms,zjmf}. We omit the details here.

\bo\lab{bo7} For $b>0$ small, there exist
$\{u_n\}_n\subset I_b^{D_b}\cap S^d$ such that
$I_b'(u_n)\rg 0$ as $n\rg\iy$.
\eo

\vskip2pt
\noindent
{\bf The completion of Proof of Theorem \ref{Theorem 1} }
\bp It follows from Proposition \ref{bo7} that there exists
$b_0>0$ such that for $b\in(0,b_0)$, there exists $\{u_n\}\in
I_b^{D_b}\cap S^d$ with $I_b'(u_n)\rg 0$ as
$n\rg\iy$. Thus, there exists $u_b\in H_r^s(\RN)$, up to a subsequence, such that $u_n\rg u_b$ weakly in  $H_r^s(\RN)$, $u_n\rg u_b$  strongly in  $L^p(\RN)$, $p\in(2,2_s^*)$ and $u_n\rg u_b$  a.e  in  $\RN$.
Next, we claim that $I_b'(u_b)=0$ for $b$ small. Set $f(t)=g(t)+t^{2_s^*-1}$. By Lemma \ref{l2}, we have
$\int_{\RN}g(u_n)\varphi=\int_{\RN}g(u_b)\varphi+o_n(1)$ for any $\varphi\in C_0^\infty (\RN)$
and $\int_{\RN}g(u_n)u_n=\int_{\RN}g(u_b)u_b+o_n(1).$
Let $v_n=u_n-u_b$ and $\|v_n\|_{\mathcal{D}^{s,2}}^2\rg A\ge 0$, then $\|u_n\|_{\mathcal{D}^{s,2}}^2=\|u_b\|_{\mathcal{D}^{s,2}}^2+A+o_n(1)$.
From $I_b'(u_n)\rg 0$, we have
\be\lab{guji1}
\left(a+b\|u_b\|_{\mathcal{D}^{s,2}}^2+bA\right)\|u_b\|_{\mathcal{D}^{s,2}}^2+\|u_b\|_2^2 =\int_{\RN}g(u_b)u_b+\|u_b\|_{2_s^*}^{2_s^*}.
\ee
The corresponding Pohoz\v{a}ev identity is
\be\lab{pohozaev3}
\frac{N-2s}2\left(a+b\|u_b\|_{\mathcal{D}^{s,2}}^2+bA\right)\|u_b\|_{\mathcal{D}^{s,2}}^2+\frac{N}2\|u_b\|_2^2=N\int_{\RN}G(u_b)+\frac{N}{2_s^*}\|u_b\|_{2_s^*}^{2_s^*}.
\ee
It follows from $I_b'(u_n)u_n\rg 0$ and Brezis-Lieb Lemma that
$$
\left(a+b\|u_b\|_{\mathcal{D}^{s,2}}^2+bA\right)(\|u_b\|_{\mathcal{D}^{s,2}}^2+A)+(\|u_b\|_2^2+\|v_n\|_2^2) =\int_{\RN}g(u_b)u_b+\|u_b\|_{2_s^*}^{2_s^*}+\|v_n\|_{2_s^*}^{2_s^*}+o_n(1).
$$
Together with (\ref{guji1}), we have
\be\lab{guji3}
\left(a+b\|u_b\|_{\mathcal{D}^{s,2}}^2+bA\right)A+\|v_n\|_2^2=\|v_n\|_{2_s^*}^{2_s^*}+o_n(1).
\ee
It follows from Lemma \ref{l2} that
$
A\le \frac{1}a\left(\frac{A}{S_s}\right)^{\frac{2_s^*}2}+o(1).
$
If $A=0$, we have done. If $A>0$, then
$A\ge a^{\frac{N-2s}{2s}}S_s^{\frac{N}{2s}}.$
%\be\lab{guji4}
%A\ge a^{\frac{N-2s}{2s}}S_s^{\frac{N}{2s}}.
%\ee
By the Pohoz\v{a}ev identity (\ref{pohozaev3}) and (\ref{guji3}),
\begin{align*}
I_b(u_n)&=\left(\frac{1}2-\frac{1}{2_s^*}\right)a(\|u_b\|_{\mathcal{D}^{s,2}}^2+A)+\left(\frac{1}4-\frac{1}{2_s^*}\right)b(\|u_b\|_{\mathcal{D}^{s,2}}^2+A)^2+\left(\frac{1}2-\frac{1}{2_s^*}\right)\|v_n\|_2^2+o(1)\\
&\ge \left(\frac{1}2-\frac{1}{2_s^*}\right)aA+b\left(\frac{1}4-\frac{1}{2_s^*}\right)(\|u_b\|_{\mathcal{D}^{s,2}}^2+A)^2+o(1).
\end{align*}
On the other hand, from $\{u_n\}\in S^d$, for $d$ small, there exist $\ti u_n\in S_r$ and $\ti v_n\in H^s(\RN)$ such that $u_n=\ti u_n+\ti v_n$ with $\|\ti v_n\|\le d$. Thus $\|u_n\|_{\mathcal{D}^{s,2}}^2\le \|\ti v_n\|_{\mathcal{D}^{s,2}}^2+\|\ti u_n\|_{\mathcal{D}^{s,2}}^2
\le 1+\sup_{v\in S_r}\|v\|_{\mathcal{D}^{s,2}}^2\triangleq B$
%\begin{align*}
%\|u_n\|_{\mathcal{D}^{s,2}}^2&\le \|\ti v_n\|_{\mathcal{D}^{s,2}}^2+\|\ti u_n\|_{\mathcal{D}^{s,2}}^2\\
%&\le 1+\sup_{v\in S_r}\|v\|_{\mathcal{D}^{s,2}}^2\triangleq B,
%\end{align*}
which implies that $\|u_b\|_{\mathcal{D}^{s,2}}^2+A\leq 2B$, where $B$ is independent of $b,n$ and $d$. So
\begin{align*}
I_b(u_n)\ge \left(\frac{1}2-\frac{1}{2_s^*}\right)aA-4b\left|\frac{1}4-\frac{1}{2_s^*}\right|B^2+o(1).
\end{align*}
Meanwhile, from $\limsup_{n\rg \infty}I_b(u_n)\le D_b$, we get
\begin{align*}
\left(\frac{1}2-\frac{1}{2_s^*}\right)aA\le D_b+b\left|\frac{1}4-\frac{1}{2_s^*}\right|B^2.
\end{align*}
Together with $A\ge a^{\frac{N-2s}{2s}}S_s^{\frac{N}{2s}}$, we have $\frac{s}N(aS_s)^{\frac{N}{2s}}\le D_b+b\left|\frac{1}4-\frac{1}{2_s^*}\right|B^2\rg E, \mbox{ as } b\rg 0,$
%\begin{align*}
%\frac{s}N(aS_s)^{\frac{N}{2s}}\le D_b+b\left|\frac{1}4-\frac{1}{2_s^*}\right|B^4\rg E, \mbox{ as } b\rg 0,
%\end{align*}
which is a contradiction with $E<\frac{s}N(aS_s)^{\frac{N}{2s}}$.
So, the claim is true. Since $u_n\in S^d$, then for $d$ small, $u_b\not\equiv 0$.
Thus, for $b$ and $d$ small, there exists $u_b\in H_r^s(\RN)$ which is a nontrivial solution of (\ref{fractional kirchhoff}). In the following, we investigate the asymptotic behavior of $u_b$ as $b\rg 0$. Noting that $D_b\rg E$ as $b\rg 0$, the similar proof as that in Proposition \ref{bo5}, we obtain that there exist $u\not\equiv 0$ such that $u_b\rg u$ strongly in $H_r^s(\RN)$ with $L'(u)=0$ and $L(u)=E$. The proof is finished.
\ep

%%%%%%%%%%%%%%%%%%%%%%%%%%%%%%%%%%%%%%%%%%%%%%%%%%%%%%%%%%%%%%%%%%%%%%%%%%%%%%%%%%%%%%%%%%%%%%%%%

\noindent{\bf Acknowledgements.}\,\,
{\rm This work is supported  by the National Natural Science Foundation of China (11271364).}

\bigskip
\bigskip

\begin{thebibliography}{10}

%\bibitem{G.Alberti}
%G. Alberti, G. Bouchitt\'{e}, P. Seppecher,
%Phase transition with the line-tenstion effect,
%{\it Arch. Ration. Mech. Anal.} {\bf 144} (1998), 1-46.

\bibitem{App} D. Applebaum,
L\'{e}vy processes-from probability theory to finance and quantum groups,
{\it Notices of the American Math Soc.},
{\bf 51} (2004), 1320-1331.


\bibitem{Alves} C. O. \ Alves, M. A. S.\ Souto, M.\ Montenegro,
Existence of a ground state solution for a nonlinear scalar field equation with critical growth,
{\it Calc. Var. PDE.} {\bf 43} (2012), 537--554.


\bibitem{Ambrosio} V. Ambrosio, T. Isernia, A multiplicity result for a fractional Kirchhoff equation in $\RN$ with a general nonlinearity.
{\it arXiv:1606.05845}

\bibitem{AutuoriFisPucci} G. Autuori, A. Fiscella and P. Pucci,
Stationary Kirchhoff problems invoving a fractional elliptic operator and a critical nonlinearity,
{\it Nonlinear Anal.}
{\bf 125}(2015), 699-714.

\bibitem{Autuori} G. Autuori and P. Pucci,
Elliptic problems involving the fractional Laplacian in $\RN$,
{\it J. Differential Equations.}
{\bf 255} (2013), 2340-2362.


%\bibitem{B.Barriosa}B. Barriosa, E. Coloradoc, R. Servadeid, F. Soriaa, A critical fractional equation with concave-convex power nonlinearities, {\it Ann. Inst. H. Poincare Anal. Nonlineaire}, 32, 875-900(2015).
%
%\bibitem{Bernstein} S. Bernstein,
%Sur une classe d'\'{e}quations fonctionnelles aux d\'{e}iv\'{e}s partielles.
%{\it Bull. Acad. Sci. URSS. S\'{e}. Math.}
%{\bf4} (1940), 17-26.

%
%\bibitem{bisci} G. M. Bisci, V. D. Radulescu, Ground state solutions of scalar field fractional Schr\"{o}dinger equation, {\it Calc. Var. Partial Differential Equatioans}, 54, 2985-3008(2015).


\bibitem{c.brandle}  C. Br\"{a}ndle,  E. Colorado, U. S\'{a}nchez, A concave-convex elliptic problem involving the fractionnal Laplacian, {\it Proc. R. Soc. Edinb.}, {\bf 143A}(2013), 39-71.

%\bibitem{BybeonJean} J. Byeon and L. Jeanjean,
%Standing waves for nonlineat Schr\"{o}dinger equations with a general nonlinearity,
%{\it Arch. Ration. Mech. Anal.}
%{\bf 185}(2007), 185-200.

%\bibitem{Byeon-Zhang-Zou}
%J. Byeon, J. Zhang, W. Zou, Singularly perturbed nonlinear Dirichlet
%problems involving critical growth,
%{\it Calc. Var. PDE.} {\bf 47} (2013), 65--85.


%\bibitem{L.A.Caffarelli}
%L. A. Caffarelli, S. Salsa, L. Silvestre, Regularity estimates for the solution and the free boundary of the obstacle problem for the fractional Laplacian,
%{\it Invent. Math.}  {\bf 171(2)} (2008),425-461.


\bibitem{Changxiaojun1} X. Chang and Z. Q. Wang,
Ground state of scalar field equations involving a fractional Laplacian with general nonlinearity, {\it Nonlinearity}, {\bf 26} (2013), 479--494.
%
%\bibitem{Changxiaojun2} X. Chang and Z. Q. Wang,
%Nodal and multiple solutions of nonlinear problems involving the fractional Laplacian,
%{\it J. Differential Equations.}
%{\bf 256}(2014), 2965-2992.
%
%\bibitem{M.chipot} M. Chipot, B. Lovat, Some remarks on nonlocal elliptic and parabolic problems,
%{\it Nonlinear Anal.}
%{\bf 30}(1997), 4619-4627.


\bibitem{Cots} A. Cotsiolis and N. K. Tavoularis, Best constants for Sobolev inequalities for higher
order fractional derivatives, {\it J. Math. Anal. Appl.,} {\bf 295} (2004), 225--236.


\bibitem{dyps} Y. Deng, S. Peng, W. Shuai,
Existence and asymptotic behavior of nodal solutions for the Kirchhoff-type problems in $\mathbb{R}^3$.
{\it J. Funct. Anal.}
{\bf 269} (2015), 3500-3527.


\bibitem{DPV}
E.\ Di Nezza, G.\ Palatucci, E.\ Valdinoci,
Hitchhiker's guide to the fractional Sobolev spaces,
\emph{Bull.\ Sci.\ Math.} \textbf{136} (2012), 512--573.
%
%\bibitem{Dip} S. Dipierro, G. Palatucci and E. Valdinoci, Existence and symmetry results for a Schr\"{o}dinger type problem
%involving the fractional Laplacian,  {\it Le Matematiche}, {\bf LXVIII} (2013), 201--216.
%
%
%\bibitem{Fiscella1} A. Fiscella,
%A fractional Kirchhoff problem involving a singular term and a critical nonlinearity.
%{\it arMiv: 1703.07861v1.}

\bibitem{Fiscella} A. Fiscella and E. Valdinoci,
A critical Kirchhoff type problem involving a nonlocal operator.
{\it Nonlinear Anal.}
{\bf 94}(2014), 156-170.

%
%
%\bibitem{hylg2} Y. He, G. Li, S. Peng,
%Concentrating Bound States for Kirchhoff type problems in $\mathbb{R}^3$ involving critical Sobolev exponents.
%{\it Adv. Nonlinear Stud.}
%{\bf 14} 2014, 441-468.

\bibitem{X.H.W.Z1} X. He, W. Zou,
Existence and concetration behavior of positive solutions for a Kirchhoff equation in $\R^3$,
{\it J. Differential Equations}
{\bf 252} (2012), 1813-1834.


\bibitem{Kirchhoff} G. Kirchhoff, Mechanik, Teubner, Leipzig, 1883.

\bibitem{Lions1} P. L. Lions, Sym\'{e}trie et compacit\'{e} dans les espaces de Sobolev, {\it J. Funct. Analysis}, {\bf 49} (1982), 315--334.

\bibitem{Lask2}
N. Laskin, {\it Fractional quantum mechanics and Levy path integrals}, Physics Letters A {\bf 268} (2000), 298--305.


\bibitem{liuzhisumarco} Z. S. Liu, M. Squassina, J. J. Zhang, Ground state for fractional Kirchhoff equations with critical nonlinearity in low Dimension. {\it arXiv:1612.07914.}


\bibitem{Nyamoradi} N. Nyamoradi,
Existence of three solutions for Kirchhoff nonlocal operators of elliptic type,
{\it Math. Commun.}
{\bf 18}(2013), 489-502.



%\bibitem{Pohokirc} S.I. Poho\u{z}aev,
%A certain class of quasilinear hyperbolic equations.
%{\it Mat. Sb.}
%{\bf 96}(1975), 152-166.


\bibitem{Pucciandsaldi} P. Pucci and S. Saldi,
Critical stationary Kirchhoff equations in $\RN$ involving nonlocal operators,
{\it Rev. Mat. Iberoam.}
{\bf 32}(2016), 1-22.

%
%\bibitem{PatriziaMingqiXiangzbl} P. Pucci, M. Q. Xiang and B. L. Zhang,
%Existence and Multiolicity of entire solutions for fractional $p$-Kirchhoff equations,
%{\it Adv. Nonlinear Anal.}
%{\bf 5}(2016), 27-55.

\bibitem{Y.Sire} Y. Sire, E. Valdinoci,
Fractional Laplacian phase transtion and boundary reactions: a geometric inequality and a symmetry result.
{\it J. Funct. Anal.} {\bf 256(6)} (2009), 1842-1864.


%\bibitem{M.Struwe} M. Struwe,
%Variational methods. Application to nonlinear partial differential equations and hamiltonian Systems,
%    {\it Springer-Verlag,}
%    (1990).
%
%\bibitem{wjtl} J. Wang, L. Tian, J. Xu, F. Zhang,
%Multiplicity and concentration of positive solutions for a Kirchhoff type problem with critical growth.
%{\it J. Differential Equatiaons.}
%{\bf 253}(2012), 2314-2351.

\bibitem{xiangzhangyang} M. Q. Xiang, B. L. Zhang and M. M. Yang, A fractional Kirchhoff-type problem in $\RN$ without the (AR) condition.
{\it Complex Var. Elliptic Equ.}
{\bf 61(11)}(2016),1481-1493.

\bibitem{zjjjmms} J. J. Zhang, J. M. do O and M. Squassina, Schr\"{o}dinger-Poisson systems with a general critical nonlinearity.
{\it Commun. Contemp. Math.},1650028 (2016).

\bibitem{zjmf} J. J. Zhang, J. M. do O and M. Squassina,
Fractional Schr\"{o}dinger-Poisson systems with a general subcritical or critical nonlinearity,
{\it Adv. Nonlinear Stud.}
{\bf 16(1)}(2016), 15-30.





\end{thebibliography}
\end{document}